\documentclass[twoside,a4paper,12pt,centertags,reqno]{amsart} 
\usepackage{amsmath,amssymb,verbatim,vmargin}
\usepackage{color}
\usepackage{tikz}
\usepackage{hyperref}
\hypersetup{colorlinks,bookmarks=true,linktocpage=true,citecolor=blue,linkcolor=magenta}

\usepackage{eulervm}   

\allowdisplaybreaks 
\usepackage{color}
\usepackage[all]{xy} 

\usepackage{mathtools}
\usepackage{MnSymbol} 

\theoremstyle{plain}
\newtheorem{thm}{Theorem}[section]

\newtheorem{con}[thm]{Conjecture}

\theoremstyle{definition}
\newtheorem{defn}[thm]{Definition}

\usepackage{enumerate}

\newcommand{\cB}{{\mathcal B}}

\newcommand{\cH}{{\mathcal H}}

\newcommand{\cK}{{\mathcal K}}
\newcommand{\cL}{{\mathcal L}}

\def\barint_#1{\mathchoice
            {\mathop{\vrule width 6pt
height 3 pt depth -2.5pt
                    \kern -9.5pt
\intop \kern -4pt}\nolimits_{#1}}%
            {\mathop{\vrule width 5pt height
3 pt depth -2.6pt
                    \kern -6.5pt
\intop \kern -4pt}\nolimits_{#1}}%
            {\mathop{\vrule width 5pt height
3 pt depth -2.6pt
                    \kern -6pt
\intop \kern -4pt}\nolimits_{#1}}%
            {\mathop{\vrule width 5pt height
3 pt depth -2.6pt
          \kern -6pt \intop \kern -4pt}\nolimits_{#1}}}
          
           \def\bariint_#1{\mathchoice
            {\mathop{\vrule width 15pt
height 3 pt depth -2.5pt
                    \kern -15.8pt
\intop \kern -8pt\intop \kern -4pt}\nolimits_{#1}}%
            {\mathop{\vrule width 9pt height
3 pt depth -2.6pt
                    \kern -10.5pt
\intop \kern -8pt\intop \kern -4pt}\nolimits_{#1}}%
            {\mathop{\vrule width 9pt height
3 pt depth -2.6pt
                    \kern -10pt
\intop \kern -8pt\intop \kern -4pt}\nolimits_{#1}}%
            {\mathop{\vrule width 9pt height
3 pt depth -2.6pt
          \kern -8pt \intop \kern -10pt\intop \kern -4pt}
      \nolimits_{  #1}}}

\def\barintlim_#1{\mathchoice
            {\mathop{\vrule width 6pt
height 3 pt depth -2.5pt
                    \kern -8.8pt
\intop \kern -4pt}\limits_{#1}}%
            {\mathop{\vrule width 5pt height
3 pt depth -2.6pt
                    \kern -6.5pt
\intop \kern -4pt}\limits_{#1}}%
            {\mathop{\vrule width 5pt height
3 pt depth -2.6pt
                    \kern -6pt
\intop \kern -4pt}\limits_{#1}}%
            {\mathop{\vrule width 5pt height
3 pt depth -2.6pt
          \kern -6pt \intop \kern -4pt}\limits_{#1}}}
          
           \def\bariintlim_#1{\mathchoice
            {\mathop{\vrule width 15pt
height 3 pt depth -2.5pt
                    \kern -15.8pt
\intop \kern -8pt\intop \kern -4pt}\limits_{#1}}%
            {\mathop{\vrule width 9pt height
3 pt depth -2.6pt
                    \kern -10.5pt
\intop \kern -8pt\intop \kern -4pt}\limits_{#1}}%
            {\mathop{\vrule width 9pt height
3 pt depth -2.6pt
                    \kern -10pt
\intop \kern -8pt\intop \kern -4pt}\limits_{#1}}%
            {\mathop{\vrule width 9pt height
3 pt depth -2.6pt
          \kern -8pt \intop \kern -10pt\intop \kern -4pt}
      \limits_{  #1}}}
          
\renewcommand{\iint}{\int \kern -3pt\int}       

\numberwithin{equation}{section}
\usepackage[symbol]{footmisc}

\title[On noncommutative H\"older for weak Schatten]{On noncommutative H\"older inequality of Sukochev and Zanin for weak Schatten class}

\author{Yi C. Huang} 
\address{School of Mathematical Sciences, Nanjing Normal University, Nanjing 210023, People's Republic of China}
\email{Yi.Huang.Analysis@gmail.com}
\urladdr{https://orcid.org/0000-0002-1297-7674}

\author{Sijie Luo}
\address{School of Mathematics and Statistics, Central South University, HNP-LAMA, Changsha 410075, People's Republic of China}
\email{sijieluo@csu.edu.cn}

\date{\today} 

\thanks{Research of the authors is partially supported by the National NSF grants of China (nos. 11801274 and 12201646)
and the Natural Science Foundation of Hunan (no. 2023JJ40696).
YCH would like to thank Shangquan Bu for kind support and Rapha\"el Ponge for helpful communications.}

\begin{document}

\begin{abstract}
Sukochev and Zanin resolved an open problem due to B. Simon 
concerning optimal constants in H\"older inequality for the weak Schatten classes of compact operators.
In this note we observe that these constants, by introducing the modified weak Schatten quasi-norms, 
can be renormalised so that the original Simon's conjecture (with optimal constant 1) does hold.
We also provide an unexpectedly simple proof for the modified H\"older inequality and its sharpness.
\end{abstract}

\maketitle


\section{Introduction}

Let $\cH$ be a complex separable Hilbert space.
Denote by $\cB(\cH)$ the set of all bounded operators on $\cH$ equipped with uniform operator norm $\|\cdot\|_\infty$
and by $\cK(\cH)$ the ideal of all compact operators on $\cH$.
For $T\in\cK(\cH)$, let $\mu(T)=\{\mu(k,T)\}_{k=0}^\infty$ be the sequence of eigenvalues of $|T|=\sqrt{T^*T}$ arranged in non-increasing order with multiplicities.
It is convenient to replace the sequence $\mu(T)$ by the function 
$$\mu(t,T)=\sum_{k=0}^\infty\mu(k,T)\mathbf{1}_{[k,k+1)}(t),\quad t\geq0,$$
and one has the following crucial property 
\begin{equation} \label{eqn:property}
\mu(t_1+t_2,T_1T_2)\leq\mu(t_1,T_1)\mu(t_2,T_2),
\end{equation}
where $T_1, T_2\in \cK(\cH)$ and $t_1, t_2\geq0$.

\begin{defn}
Define for $0<p<\infty$ the weak Schatten class $\cL_{p,\infty}$ by those $T$ with
$$\|T\|_{p,\infty}:=\sup_{t>0}\,t^{\frac1p}\mu(t,T)<\infty.$$
\end{defn}

For further materials about the classes $\cL_{p,\infty}$ (in particular $\cL_{1,\infty}$) and their applications, 
see Connes \cite{Con94}, Lord-Sukochev-Zanin \cite{LorSukZan13} and Ponge \cite{Pon23}.

B. Simon asked (see \cite[p. 32]{Sim79} and \cite[p. 21]{Sim05}) the following

\begin{con}
Let $p, q, r>0$ be such that $\frac1r=\frac1p+\frac1q$. Is it true that
$$\|TS\|_{r,\infty}\leq \|T\|_{p,\infty}\|S\|_{q,\infty},\quad T\in\cL_{p,\infty}, S\in\cL_{q,\infty}?$$
\end{con}

This was answered in the negative by Sukochev and Zanin in \cite{SukZan21}. They proved
\begin{thm}[Sukochev-Zanin]
Let $p, q, r>0$ be such that $\frac1r=\frac1p+\frac1q$. Then
$$\|TS\|_{r,\infty}\leq \frac{(p+q)^{\frac1p+\frac1q}}{q^{\frac1p}p^{\frac1q}}\|T\|_{p,\infty}\|S\|_{q,\infty},\quad T\in\cL_{p,\infty}, S\in\cL_{q,\infty},$$
and the constant is optimal.
\end{thm}

For other versions of noncommutative H\"older inequalities, 
see Guido-Isola \cite{GuiIso03}, Dykema-Skripka \cite{DykSkr15} and Sukochev \cite{Suk16}.
The aim of this note is to point out the following amusing observation: 
the Sukochev-Zanin constants can be written as
$$\frac{(p+q)^{\frac1p+\frac1q}}{q^{\frac1p}p^{\frac1q}}=\frac{p^{\frac1p}q^{\frac1q}}{r^{\frac1r}}.$$
This motivates us to introduce for $T\in\cL_{p,\infty}$ the following renormalised quantity 
\begin{equation} \label{eqn:renorm}
\|T\|'_{p,\infty}:=\sup_{t>0}\,(pt)^{\frac1p}\mu(t,T)=\sup_{t>0}\,t^{\frac1p}\mu\left(\frac{t}{p},T\right).
\end{equation}
So the original conjecture of Simon does hold for the scale of quantities $\|\cdot\|'_{p,\infty}$.

\begin{thm}[$\simeq$ Sukochev-Zanin] \label{thm:FZ}
Let $p, q, r>0$ be such that $\frac1r=\frac1p+\frac1q$. Then
$$\|TS\|'_{r,\infty}\leq \|T\|'_{p,\infty}\|S\|'_{q,\infty},\quad T\in\cL_{p,\infty}, S\in\cL_{q,\infty},$$
and the constant is optimal.
\end{thm}

Our proof for this result is unexpectedly simple
and we believe that the renormalised quantities \eqref{eqn:renorm} would be useful in other quantitative studies of $\cL_{p,\infty}$.

\section{Proof of Theorem \ref{thm:FZ}}

Let $p, q, r>0$ be such that $\frac1r=\frac1p+\frac1q$. 
Let $T\in\cL_{p,\infty}$ and $S\in\cL_{q,\infty}$. Thus, by \eqref{eqn:property},
$$t^{\frac1r}\mu\left(\frac{t}{r},TS\right)=t^{\frac1p+\frac1q}\mu\left(\frac tp+\frac tq,TS\right)\leq t^{\frac1p}\mu\left(\frac{t}{p},T\right)\,\, t^{\frac1q}\mu\left(\frac{t}{q},S\right).$$
Taking the supremum for $t>0$ completes the proof of the H\"older inequality.

For the sharpness of the H\"older inequality, it suffices to note that in above reasoning, the inequality is saturated by the particular choices
where $T$ (resp. $S$) is given by the rank $n$ diagonal operator with diagonal values $\{k^{-\frac1p}\}_{1\leq k\leq n}$ (resp. $\{k^{-\frac1q}\}_{1\leq k\leq n}$).

\bigskip

\section*{\textbf{Compliance with ethical standards}}

\bigskip

\textbf{Conflict of interest} The authors have no known competing financial interests
or personal relationships that could have appeared to influence this reported work.

\bigskip

\textbf{Availability of data and material} Not applicable.

\bigskip

\bibliographystyle{alpha}

\bibliography{HuaY-LuoSJ-HolderSchatten} 

\begin{thebibliography}{LSZ13}

\bibitem[Con94]{Con94}
Alain Connes.
\newblock {\em Noncommutative geometry}.
\newblock Academic Press, Inc., San Diego, CA, 1994.

\bibitem[DS15]{DykSkr15}
Ken Dykema and Anna Skripka.
\newblock H\"{o}lder's inequality for roots of symmetric operator spaces.
\newblock {\em Studia Math.}, 228(1):47--54, 2015.

\bibitem[GI03]{GuiIso03}
Daniele Guido and Tommaso Isola.
\newblock Dimensions and singular traces for spectral triples, with
  applications to fractals.
\newblock {\em J. Funct. Anal.}, 203(2):362--400, 2003.

\bibitem[LSZ13]{LorSukZan13}
Steven Lord, Fedor Sukochev, and Dmitriy Zanin.
\newblock {\em Singular traces}, volume~46 of {\em De Gruyter Studies in
  Mathematics}.
\newblock De Gruyter, Berlin, 2013.
\newblock Theory and applications.

\bibitem[Pon23]{Pon23}
Rapha\"{e}l Ponge.
\newblock Connes' integration and {W}eyl's laws.
\newblock {\em J. Noncommut. Geom.}, 17(2):719--767, 2023.

\bibitem[Sim79]{Sim79}
Barry Simon.
\newblock {\em Trace ideals and their applications}, volume~35 of {\em London
  Mathematical Society Lecture Note Series}.
\newblock Cambridge University Press, Cambridge-New York, 1979.

\bibitem[Sim05]{Sim05}
Barry Simon.
\newblock {\em Trace ideals and their applications}, volume 120 of {\em
  Mathematical Surveys and Monographs}.
\newblock American Mathematical Society, Providence, RI, second edition, 2005.

\bibitem[Suk16]{Suk16}
Fedor Sukochev.
\newblock H\"{o}lder inequality for symmetric operator spaces and trace
  property of {K}-cycles.
\newblock {\em Bull. Lond. Math. Soc.}, 48(4):637--647, 2016.

\bibitem[SZ21]{SukZan21}
Fedor Sukochev and Dmitriy Zanin.
\newblock Optimal constants in non-commutative {H}\"{o}lder inequality for
  quasi-norms.
\newblock {\em Proc. Amer. Math. Soc.}, 149(9):3813--3817, 2021.

\end{thebibliography}
 
\end{document}